\documentclass[a4paper,10pt]{amsart}
\usepackage[usenames]{color} 
\usepackage{amssymb} 
\usepackage{amsmath} 
\usepackage{amsthm}
\usepackage{mathrsfs}
\usepackage[all]{xy}
\usepackage[utf8]{inputenc} 
\newtheorem{theo}{Theorem}[section]

\newtheorem{exam}{Examples}[section]
\newtheorem{prop}{Proposition}[section]
\newtheorem{cor}{Corollary}[section]
\newtheorem{rem}{Remark}[section]
\theoremstyle {remark}

\title{Smooth crossed product of minimal unique ergodic diffeomorphism of odd sphere}
\author{Liu, Hongzhi}

\address{Department of Mathematics\\Jilin University\\Changchun 130012\\P.R. CHINA}

\email{hzliu13@mails.jlu.edu.cn}

\thanks{Research, supported by NNSF of China (11201171), NNSF of China (11531003), NNSF of China (11271150) and NNSF of China (11401088)}

\begin{document}

\pagestyle{plain}
\begin{abstract}
For minimal unique ergodic diffeomorphisms $\alpha_n$ of $S^{2n+1} (n>0)$ and $\alpha_m$ of $S^{2m+1}(m>0)$, the $C^*$-crossed product algebra $C(S^{2n+1})\rtimes_{\alpha_n} \mathbb{Z}$ is isomorphic to $C(S^{2m+1})\rtimes_{\alpha_m} \mathbb{Z}$ even though $n\neq m$ . However, by cyclic cohomology, we show that smooth crossed product algebra $C^\infty(S^{2n+1})\rtimes_{\alpha_n} \mathbb{Z}$ is not isomorphic to $C^\infty(S^{2m+1})\rtimes_{\alpha_m} \mathbb{Z}$ if $n\neq m$. 
\end{abstract} 


\date{\today}
\maketitle 
\section{Introduction}

$C^*$-algebra classification theory can be used to study dynamical systems. Pimsner, Voiculescu (\cite{pv}) and Rieffel (\cite{rie}) proved that two irrational rotation dynamical systems are flip conjugate to each other if and only if their corresponding irrational rotation $C^*$-algebras are isomorphic to each other. Giordano, Putnam and Skau have shown that the minimal dynamical systems of Cantor set can be completely classified by $C^*$-crossed product algebras up to strong orbit equivalence (\cite{gps}). See \cite{lint}, \cite{Hlin}, \cite{ch3}, \cite{ch4}, \cite{win} for more examples.

However, there are examples of different minimal diffeomorphisms give the same $C^*$-algebras. Let $\alpha_l$ be minimal unique ergodic diffeomorphism of $S^{2l+1}, l=1,2,\dots $. It is well known that the ordered $K$-theory of $C(S^{2n+1})\rtimes_{\alpha_n} \mathbb{Z}$ and $C(S^{2m+1})\rtimes_{\alpha_m} \mathbb{Z}$ are isomorphic to each other (\cite{ch5}). This implies 
\[C(S^{2n+1})\rtimes_{\alpha_n} \mathbb{Z}\cong C(S^{2m+1})\rtimes_{\alpha_m} \mathbb{Z}\] no matter if $n=m$ or not according to a theory of Toms and Winter (\cite{win}) and Phillips (\cite{ch5}). See \cite{ch1} and \cite{ch2} for more examples.

Let $M$ be a compact manifold, Chern map naturally defines a graded structure for topological $K$-theory of $M$: 
\[K^0(M) \stackrel{Chern}{\to } H^0(M)\oplus H^2(M)\dots \oplus H^{2n}(M)\dots,\]
\[K^1(M) \stackrel{Chern}{\to } H^1(M)\oplus H^3(M)\dots \oplus H^{2n+1}(M)\dots.\]

Unfortunately, this is not the case for $K$-theory of $C^*$-algebras, where an order structure is the best one can get. Consider the example constructed by Goodearl  (\cite{good}). The classification result is obtained by Elliott and Gong (\cite{EG2}).
\begin{exam}
	Let $M$ be a connected finite dimensional compact manifold, $\{x_i\}_{i=1}^\infty$ be a dense subset of $M$. Define an inductive limit system of $C^*$-algebra as follows
	\[C(M)\stackrel{\Phi_0}{\longrightarrow}M_2(C(M))\stackrel{\Phi_1}{\longrightarrow}M_4(C(M))\dots \stackrel{\Phi_{n-1}}{\longrightarrow}M_{2^n}(C(M))\stackrel{\Phi_{n}}{\longrightarrow}\dots\]
	
	where 
	\[\Phi_i(f)(x)= \left(\begin{array}{cc}
	                f(x)&                         \\
	                 &f(x_i)                   \\
	                \end{array}
	                \right). \]
	                
	 It is proved by Goodearl that such limit $C^*$-algebras are of real rank zero. Obviously two different choices of dense sets $\{x_i\}_{i=1}^{\infty}$ give two shape equivalent inductive limit systems and therefore give a same $C^*$-algebra by Theorem 2.2 of \cite{good}, that is, the limit $C^*$-algebra does not depend on the choice of the dense subset $\{x_i\}_{i=1}^{\infty}$. We denote the limit $C^*$-algebras by $\mathbb{A}(M)$. Take $M$ to be $S^1$ and $S^3$. In every finite stage, the nontrivial odd $K$-elements reflect different levels of cohomology of spaces, namely the Chern Character of the elements live in $H^1(S^1)$ or $H^3(S^3)$. However the Elliott invariants are the same for $\mathbb{A}(S^1)$ and $\mathbb{A}(S^3)$, which implies that $\mathbb{A}(S^1) \cong \mathbb{A}(S^3)$ . 
 
\end{exam}

$K$-theory is the only (co)homology theory that can be properly generalized from topology to the theory of $C^*$-algebra (non-commutative topology in Alain Connes's sense ). If we consider smooth algebras instead of $C^*$-algebras (certainly the dynamical systems would be restricted to smooth ones), cyclic cohomology and its graded structure become applicable as a new tool. Elliott and Gong (\cite{EG}) showed that not all ``continuous" homomorphism from $C(S^3)$ to irrational rotation algebras can be ``approximated" by ``smooth" ones in light of cyclic cohomology. This inspires us to look at the theory of cyclic cohomology.

Alain Connes invented cyclic cohomology in \cite{Al1}. It can be viewed as a generalization of deRham homology. This article can be used to examine the minimal unique ergodic diffeomorphisms of odd spheres and their smooth crossed product algebras. This example shall demonstrate that the ``homology theory" (cyclic cohomology) is no longer invariant for ``topological structure" ($C^*$-algebras) but depends on specific ``geometrical structure" (smooth algebras). 

\textit{Acknowledgments.} The author would like to thank professor Gong, Guihua for his advice of this question. Thanks go to professor Ryszard Nest both for his patient help in Shanghai and his work in 1980s. I am also very grateful to G. Elliott for his help.
\section{Preliminary}

Our strategy is to compute the graded structure of cyclic cohomology
$HC^*(\mathcal{A})$ of smooth crossed product algebra given by diffeomorphism of a compact manifold. One can immediately see from the definition that it is invariant under smooth flip conjugation.  In this section, we will introduce all the necessary notions. All algebras involved in this article are over $\mathbb{C}$.

\subsection{Smooth crossed product}
 
Convex topological algebra is algebra endowed with locally convex topology. Smooth crossed product algebra is convex topological algebra. 

Let $M$ be a finite dimensional compact manifold. Choose finitely many vector fields $X_1$, $X_2$, $\dots$, $X_n$ on $M$ which can span the tangent space at any point ($n$ is not necessarily equal to the dimension of $M$). Define seminorms $\| \bullet \|_n$:
 \[\|f\|_n=\sum_{1\leq k_1\leq \dots \leq k_n\leq n}\|X_{k_n}X_{k_{n-1}}\dots X_{k_1}f\|_\infty, n\in \mathbb{Z}_+\cup \{0\}, f\in C^\infty (M).\] 
Let $\alpha$ be a minimal unique ergodic diffeomorphism of $M$. Let $C^\infty (M)_\alpha[u,u^{-1}]$ be the algebraic crossed product of $C^\infty (M)$ by $\mathbb{Z}$. Let $\|\alpha^t\|_i$ be the operator seminorms defined by $\|\bullet\|_i$ on $C^\infty (M)$, i.e.  
\[\|\alpha^t\|_i \triangleq \sup_{f\in C^\infty(M), \|f\|_i =1} \|\alpha^t(f)\|_i.\]
Define a sequence of maps 
\[\rho_k:\mathbb{Z} \to  \mathbb{R}^+, k=1,2,\dots ,\]
by 
$\rho_k(n)  = \nolinebreak \sup_{i\leqslant k}(\sum_{t=-n}^n\|\alpha^t\|_i)^k.$

Endow $C^\infty (M)_\alpha[u,u^{-1}]$ with the topology defined by the following  seminorms:
\[ \|\sum_n f_n u^n\|_k=\sup_n \rho_k(n)\|f_n\|_k, f_n \in C^\infty (M).\]
This topology does not depend on the choice of $X_1$, $X_2$,  $\dots$, $X_n$. Then the completion of $C^\infty (M)_\alpha[u,u^{-1}]$ is the smooth crossed product algebra $C^\infty(M)\rtimes_\alpha \mathbb{Z}$. 

\begin{rem}
 Here we adopt the topology for smooth crossed product algebra given by Nest (\cite{nest}) instead of the one considered by Schweitzer (\cite{sch}) and Phillips (\cite{ch1}). 
 \end{rem}
 
 \begin{rem}
In fact,  $C^\infty(M)\rtimes_\alpha \mathbb{Z}$ is a Fréchet $*$-algebra.
 \end{rem}
 
\subsection{Cyclic cohomology}

Without loss of generality, we state cyclic cohomology theory for algebra (over $\mathbb{C}$) only.

Recall the definition of Hochschild $n$ cochain:
 \[C^n (\mathcal{A})= Hom(\mathcal{A}^{\otimes(n+1)}, \mathbb{C}), n=0,1,\dots .\]
Let $b: C^n(\mathcal{A}) \to C^{n+1}(\mathcal{A})$ be the Hochschild differential:
\begin{displaymath} 
\begin{array}{c}
(b\phi) (a_0, a_1, \dots, a_{n+1}) = \sum ^n _{i=0} (-1)^i \phi (a_0, a_1, \dots,a_ia_{i+1}, \dots, a_{n+1})
\\
\phantom{(b\phi) (a_0, a_1, \dots, a_{n+1}= \sum ^n _{i=0} (-1)^i} +(-1)^{n+1}\phi(a^{n+1}a^0, a^1,\dots,a^n).
\end{array}
\end{displaymath}
The Hochschild cohomology $H^*(\mathcal{A}, \mathcal{A}^*)$ is then the cohomology group of the complex $(C^*(\mathcal{A}), b)$.  

Let $C_{\lambda}^n(\mathcal{A})$ be the set of cyclic $n$-cochains, which are those elements in $C^n (\mathcal{A})$  satisfying the condition
\[\phi (a_n , a_0,\dots , a_{n-1})= (-1)^n \phi (a_0 ,a_1,\dots , a_n ).\]
The cyclic cohomology $HC^n (\mathcal{A})$ of $\mathcal{A}$ is then the cohomology group of the complex $(C_{\lambda}^n(\mathcal{A}), b)$. For example (\cite{Al1})
\[ HC^{2n}(\mathbb{C})=\mathbb{C}, HC^{2n+1}(\mathbb{C})=0, n = 0, 1, 2, \dots.\]

\begin{rem}
Note that when it comes to convex topological algebra, one should replace all the cochains by continuous ones. This shouldn't cause any confusions.	
\end{rem}

\subsection{The map $S$ and periodic cyclic cohomology}

Let us recall the notion of $n$-triple $(\Omega, d,\int)$ on 
$\mathcal{A}$ with homomorphism $\rho$. Let $\Omega=\oplus_{0}^{n}\Omega_i$ be a graded algebra, $\rho : \mathcal{A}\to \Omega_0$ be a homomorphism. Let $d$ be a graded derivation of degree 1 with $d^2=0$ and $\int: \Omega_n \to \mathbb{C}$ be a closed graded trace. Then $(\Omega, d,\int)$ is called an $n$-triple if they satisfy the following conditions, \\
$\mathnormal{(1)}$  $\Omega_i \times \Omega_j\subset \Omega_{i+j}, \forall i,j\in {0,1,2,\dots,n}, i+j\leq n$.\\
$\mathnormal{(2)}$  $d\Omega_i\subset \Omega_{i+1}, d(\omega\omega')=(d\omega)\omega'+(-1)^{deg\omega}\omega d\omega' , d^2=0$.\\
$\mathnormal{(3)}$  $\int d\omega =0, \int \omega'\omega = (-1)^{deg\omega deg\omega'}\int \omega \omega'$.

The tensor product of two triples $\rho: \mathcal{A}\to \Omega_0$, $(\Omega, d,\int)$ and $\rho':\mathcal{A'}\to \Omega_0'$, $(\Omega', d',\int')$ are given by $\rho\otimes\rho':\mathcal{A}\otimes \mathcal{A'}\to \Omega_0\otimes \Omega_0'$, $(\Omega\otimes \Omega',d\otimes d',\int\otimes\int')$. Recall that $d\otimes d'(\omega\otimes \omega')= d(\omega) \otimes \omega'+ (-1)^{deg{\omega}}\omega \otimes d'(\omega')$, $\int\otimes \int'(\omega\otimes \omega')= \int(\omega) \int'(\omega')$.

For any algebra $\mathcal{A}$ we have the following triple $(\Omega(\mathcal{A}), d, \int)$ of $n$-dimension over $\mathcal{A}$ with homomorphism $\rho$:\\
$\mathnormal{(1)}$ Adjoin a unit to $\mathcal{A}$ no matter if there has been one:
\[\mathcal{A}^{+}=\{a+\lambda I; a\in\mathcal{A},\lambda\in \mathbb{C}\}.\]
Define $\Omega(\mathcal{A})$ to be $\oplus_{0}^{\infty}\Omega_i(\mathcal{A})$ where 
\[\Omega_i(\mathcal{A})=\mathcal{A}^{+}\otimes \otimes^i\mathcal{A}.\]
 $\Omega(\mathcal{A})$ is usually called the universal graded algebra associated to $\mathcal{A}$. $\rho$ is then the natural inclusion.\\
$\mathnormal{(2)}$ Define the differential homomorphism $d$ from $\Omega_i(\mathcal{A})$ to $\Omega_{i+1}(\mathcal{A})$ as 
\[d((a^0 +\lambda^0I)\otimes a^1\otimes\dots \otimes a^n)=I\otimes a^0 \otimes \dots \otimes a^n.\] 
One can directly verify that $d^2=0$.\\
$\mathnormal{(3)}$ Define the product $\Omega_i(\mathcal{A}) \times \Omega_j(\mathcal{A})\to \Omega_{i+j}(\mathcal{A})$ as follows. There is a right $\mathcal{A}$-module structure on $\Omega(\mathcal{A})$ defined by the equation
\[(a^+\otimes a^1 \otimes \dots \otimes a^n)a=\sum_0^n(-1)^{n-j}a^+\otimes a^1 \dots\otimes a^ja^{j+1}\otimes\dots\otimes a.\]
This right action can be extended to an $\mathcal{A}^+$ action on $\Omega(\mathcal{A})$. Then the definition of $\Omega_i(\mathcal{A}) \times \Omega_j(\mathcal{A})\to \Omega_{i+j}(\mathcal{A})$ is given by
\[\omega(b^+ \otimes b^1 \otimes \dots \otimes b^j)=\omega b^+ \otimes  b^1 \otimes \dots \otimes b^j\,\,\, \forall \omega \in \Omega_i .\]
$\mathnormal{(4)}$  $\int: \Omega_n(\mathcal{A}) \to \mathbb{C}$ is a closed graded trace. Its existence is guaranteed by the following proposition.

\begin{prop}[\cite{Al1}] The following are equivalent: \\
$\mathnormal{(1)}$ 	$\tau$ is a closed cyclic $n$ cochain.\\
$\mathnormal{(2)}$  There is an $n$ triple on $\mathcal{A} $ with homomorphism $\rho$, s.t. 
	\[\tau(a^0, a^1,\dots a^n)=\int \rho(a^0) d \rho(a^1) \dots d \rho(a^n).\]
$\mathnormal{(3)}$  	There is an $n$ dimensional closed graded trace $\int$ on the universal grading algebras with the natural inclusion as the homomorphism $\rho$, s.t.
	\[\tau(a^0, a^1,\dots a^n)=\int a^0 d a^1 \dots d a^n.\]
\end{prop}

The point of the universal grading algebra is that if $(\Omega',d',\int')$
is a triple on $\mathcal{A}$ with a homomorphism $\rho'$, then there is a homomorphism from $\Omega(\mathcal{A})$ to $\Omega'$.

Now let us recall the notion of cup product. Represent two closed cyclic cochains $\phi \in Z_\lambda^n(\mathcal{A}) $, $\psi \in Z_\lambda^m(\mathcal{B})$ by triples on the universal grading algebra according to proposition 2.1. Note that there is always a natural homomorphism 
\[\pi: \Omega_{n+m}(\mathcal{A}\otimes\mathcal{B})\to \Omega_n(\mathcal{A})\otimes  \Omega_m(\mathcal{B}).\]
Then $\phi \in Z_\lambda^n(\mathcal{A}) $, $\psi \in Z_\lambda^m(\mathcal{B})$ defines an element in $Z_\lambda^n(\mathcal{A}\otimes \mathcal{B})$ as
\[\phi \cup \psi \triangleq (\phi\otimes\psi)\circ \pi.\]
As shown in \cite{Al1}, this formula actually defines a cup product on cyclic cohomology level: 
\[HC^n(\mathcal{A})\cup HC^m(\mathcal{B})\to HC^{n+m}(\mathcal{A}\otimes\mathcal{B}).\]
 
Let $\Delta$, $\Delta(1, 1, 1,)=1$, be the generator of $HC^2(\mathbb{C})$. $\Delta$ gives the next a 2 periodic homomorphism:
\begin{eqnarray*}
	S: HC^n(\mathbb{C}) & \to & HC^{n+2}(\mathbb{C})\\
	S(\phi) & = & \phi\cup \Delta
\end{eqnarray*}

Thus there would be two inductive limit systems:
\[ HC^0(\mathbb{C})\dots \to HC^{2n}(\mathbb{C}) \to  HC^{2n+2}(\mathbb{C})\to \dots,\]
\[ HC^1(\mathbb{C})\dots \to HC^{2n+1}(\mathbb{C}) \to  HC^{2n+3}(\mathbb{C})\to \dots.\]
The limit groups are the so called periodic cyclic cohomology
\[HP^i (\mathcal{A}) \triangleq \lim_\to HC^{2n+i}(\mathcal{A}),i=0,1.\]

Let $S(HC^*(\mathcal{A}))\subset HP^*(\mathcal{A})$ be the ultimate  image of $HC^*(\mathcal{A})$ in $HP^*(\mathcal{A})$.  $S(HC^n(\mathcal{A}))/S(HC^{n-2}(\mathcal{A}))$ actually defines a grading structure of $HP^*(\mathcal{A})$.
\subsection{Six term exact sequence}

Let $M$ be a compact manifold, $\alpha $ be a self-diffeomorphism of it, Nest obtained:
\begin{theo} [\cite{nest}]  
\begin{displaymath} 
 \xymatrix
 {HP^{0} (C^\infty (M)) \ar[r] &HP^{1}(C^\infty(M)) \rtimes_{\alpha}\mathbb{Z} ) \ar[r]& HP^{1}(C^\infty(M)) \ar[d]^{1-\alpha} \\
  HP^{0} (C^\infty (M)) \ar[u]^{1-\alpha}  & HP^{0}(C^\infty(M) \rtimes_{\alpha}\mathbb{Z} )\ar[l] & HP^{1}(C^\infty(M)) \ar[l] }.
 \end{displaymath}
 
\end{theo}
\subsection{Exact couple and spectral sequence}

Let $I: C_\lambda ^n(\mathcal{A})\to C^n(\mathcal{A})$ be the natural inclusion. Let $B_0$ be the homomorphism defined as:
\[(B_0\phi)(a^0,\dots,a^n)=\phi(I,a^0,\dots,a^n)-(-1)^{n+1}\phi(a^0,\dots,a^n,I),\]
$A$ is a homomorphism defined as:
\[A\phi = \sum_{\gamma\in \Gamma}\varepsilon(\gamma)\phi^\gamma,\]
where $\Gamma$ is the cyclic permutation group of $\{0, 1, 2\dots n\}$ as usual. Take $B$ to be $A\circ B_0$.

As proved in \cite{Al1}, there exists an exact couple:
\begin{displaymath}
\xymatrix {HC^*(\mathcal{A})\ar[rr]^S& & HC^*(\mathcal{A}) \ar[dl]^I\\
&H^*(\mathcal{A},\mathcal{A}^*) \ar[ul]^B&.}
 \end{displaymath}

Let $d_k$ be the differential homomorphism :
  \begin{displaymath}
\xymatrix{
 H^{n-2k}(\mathcal{A,A^*}) \ar[r]^B & HC^{n-2k-1}(\mathcal{A}) \ar[d]_S \ar[r]^I                               & H ^{n-2k-1}(\mathcal{A,A^*})  \\
     &  \vdots \ar[ur]_{d_k} \ar[d]_S                                              &                     \\
 H^n(\mathcal{A,A^*}) \ar[r]^B \ar[ur]  &HC^{n-1}(\mathcal{A}) \ar[r]^I                                                   & H^{n-1}(\mathcal{A},\mathcal{A^*}).}
 \end{displaymath}

Then the spectral sequence induced from the exact couple can be listed as  \\
$\mathnormal{1)}$ $(E_0^*(\mathcal{A}),d_0)= (H^*(\mathcal{A},\mathcal{A^*}),IB)$,  \\
$\mathnormal{2)}$ $E_n^*(\mathcal{A})$ is the homology of  $(E_{n-1}^*(\mathcal{A}),d_{n-1})$.

Consider the double complex $C^{m,n}=C^{m-n}(\mathcal{A}), \forall n,m \in \mathbb{Z}$. Define two homomorphisms  
$ \partial_1: C^{n,m}\rightarrow C^{n+1,m}$, and 
$ \partial_2: C^{n,m}\rightarrow C^{n-1,m}$ 
as follows:
\begin{eqnarray*} \partial_1(\phi)&=& (n-m+1)b\phi,\\
 \partial_2(\phi)&=& \left\{ \begin{array}{ll}
 \frac{1}{n-m}B\phi  & m\neq n\\
 0         &  m=n
\end{array}   \right. .
\end{eqnarray*}

\begin{theo}[\cite{Al1}]
$(\mathnormal{1})$ Let $F^q(C)=\sum_{m\geq q}C^{n,m}$ be the filtration according to $m$, then $H^p(F^q (C))=H_\lambda^n(\mathcal{A})$.	
\\
$(\mathnormal{2})$ The cohomology of the double complex $C$ is 
\[H^{2n}(C)=HP^0(\mathcal{A}),\]
\[H^{2n+1}(C)=HP^1(\mathcal{A}).\]
\\
$(\mathnormal{3})$The spectral sequence associated to filtration by $m$ is convergent. Actually it converges to $F^pHP^*(\mathcal{A})/F^{p+1}HP^*(\mathcal{A})$, which is $S(HC^n(\mathcal{A}))/S(HC^{n-2}(\mathcal{A}))$. And this spectral sequence coincides with the one induced from the exact couple. 
\end{theo}
From this theorem we know that the spectral sequence $E_n^*(\mathcal{A})$ is convergent. Denote the limit group as $E^n_\infty(\mathcal{A})$. 
\begin{cor}
\begin{eqnarray*}
E_\infty^n(\mathcal{A})&=&S(HC^{n}(\mathcal{A}))/S(HC^{n-2}(\mathcal{A})),\\
\oplus_{n=0}^\infty E_\infty^{2n}(\mathcal{A})&=& HP^0(\mathcal{A}).\\
\oplus_{n=0}^\infty E_\infty^{2n+1}(\mathcal{A})&=& HP^1(\mathcal{A}).
\end{eqnarray*}
\end{cor}

\begin{rem}
 $E_\infty^n(\mathcal{A})$ is also called the deRham homology of $\mathcal{A}$, for example, see \cite{nest}.
\end{rem}

\subsection{How to compute $E_\infty^n$ in our case}

Let $\alpha$ be a diffeomorphism of $M$. In this subsection we present a way to compute $E_\infty ^n(C^{\infty}(M) \rtimes _ \alpha\mathbb{Z})$ developed in \cite{nest}. Let $\Psi_n$ be the space of $n$-th deRham currents of $M$, $\partial$ be the usual boundary map, 
\[H^n_{eq}(M, \alpha) \triangleq homology\: group\: of\: (Ker(1-\alpha) |\Psi_n , \partial),\] 
\[H^n_{coeq}(M, \alpha)\triangleq homology\: group\: of\:(Coker(1-\alpha) | \Psi_n, \partial).\]
\begin{theo}[\cite{nest}]  
$E_\infty ^n(C^{\infty}(M) \rtimes_\alpha\mathbb{Z}) = H^n _{eq}(M,\alpha) \oplus H^{n-1} _{coeq}(M,\alpha)$.
\end{theo}

Denote the deRham homology as $H_n(M)$, deRham cohomology as 
$H^n(M)$. Note that $H_k(M)=H^{n-k}(M)$ for an $n$-dimensional manifold by Poincaré duality. Recall that 
\[HP^0(M)\cong \sum_{n=0}^{\infty}H_{2n}(M)\cong \sum_{n=0}^{\infty}H^{2n}(M),\]
\[HP^1(M)\cong \sum_{n=0}^{\infty}H_{2n+1}(M)\cong \sum_{n=0}^{\infty}H^{2n+1}(M).\]
as shown in \cite{Al1}.

\section{Main result}

Consider the minimal unique ergodic diffeomorphisms $\alpha_l$ (each one of them has to be orientation preserving diffeomorphism by Lefschetz fixed point theorem) of $S^{2l+1}$, $l=1,2,\dots$.  Their existence are proved by A. Fathi, M. R. Herman in \cite{fathi} and A. Windsor in \cite{wind}. $C^\infty (S^{2l+1}) \rtimes _ {\alpha_l} \mathbb{Z}$ is the smooth crossed product we have defined in section 2. 

\begin{exam}
$C^\infty (S^{2n+1}) \rtimes _ {\alpha_n} \mathbb{Z}$ is not isomorphic to $C^\infty (S^{2m+1}) \rtimes _ {\alpha_m}\mathbb{Z}$ if  $n\neq m, n,m=1,2,\dots$.
\end{exam}

\begin{theo}
 \begin{displaymath}
 E_\infty ^{2k+1}(C^\infty(S^{2n+1})\rtimes_{\alpha_n} \mathbb{Z}) = \left\{ \begin{array}{ll}
 \mathbb{C}, &k=0,n,\\
 {0}      ,   &else
\end{array}   \right..  
\end{displaymath}                                                                                                                                                                                                                   
 \end{theo}

\begin{proof}

The identity map on deRham currents descends to two homomorphisms:
 \begin{eqnarray*}
 \gamma_k &:& H_{eq}^k (S^{2n+1}, \alpha) \longrightarrow Ker(1-\alpha) | H_k (S^{2n+1}),\\
\beta _k &:& H_k (S^{2n+1})/ {(1-\alpha) H_k (S^{2n+1}) }\longrightarrow  H_{coeq}^k (S^{2n+1},\alpha).
  \end{eqnarray*}
Define homomorphisms $s_1$, $s_2$, $s_3$ as follows: \\
$\mathnormal{a})$ Given $\psi \in Ker(1-\alpha)|H_k (S^{2n+1})$,  there exists a $\phi$ such that  $(1-\alpha )\psi=\partial \phi $, since $(1-\alpha) \psi$ is a boundary. Define 
 \[s_1 : Coker \gamma _{k-1}\to  Coker\beta_k ,\]
 by	$s_1(\psi )=[\phi]$. \\
$\mathnormal{b})$  Given $\psi \in H_{coeq}^n(S^{2n+1},\alpha)$, there exists $\phi$ such that $\partial \psi = (1-\alpha)\phi$. Define  
\[ s_2 :Coker\beta_k  \to Ker\gamma_{k-2} ,\]
  by     $s_2(\psi) =[\partial \phi]$.
 \\
$\mathnormal{c})$  Given $\psi \in Ker\gamma_{k-2}$, then $\partial \psi = (1-\alpha) \psi = 0$ and $\psi=\partial \phi$ for some current $\phi$.  $(1-\alpha)\phi$ is closed since $\partial(1-\alpha)\phi=(1-\alpha)\partial\phi=0$. Define  
\[ 	s_3 :Ker\gamma_{k-2} \to  Ker\beta_{k-1},\]
       $s_3(\psi)  =(1-\alpha) [\phi]$.
 
It is easy to see that they are all well defined. Using nothing but basic computation, one can verify that the following sequence is exact:
\[ 0 \longrightarrow  Coker\gamma _{k-1}
 \stackrel{s_1}{\longrightarrow} Coker\beta_k\stackrel{s_2}{\longrightarrow} Ker\gamma_{k-2} \stackrel{s_3}{\longrightarrow} Ker\beta_{k-1} \longrightarrow 0 .\]
$\mathnormal{1)}$ $\Psi_k$ are automatically $\{0\}$ when $k>2n+1$, so are the
$H_{eq}^k (S^{2n+1},\alpha)$ and  $H_{coeq}^k (S^{2n+1}, \alpha)$.
This implies
 \begin{eqnarray*}
 Coker\beta_k &=&\{0\},k>2n+1,\\
 Ker\gamma_k &=&\{0\},k>2n+1.  
 \end{eqnarray*} 

 Applying the exact sequence above, we have
 \begin{eqnarray*} 
 Ker\gamma_{2n+1}&=&\{0\},\\
 Coker\gamma _{2n+1}&=&\{0\}.
 \end{eqnarray*}
 i. e. 
  \[H^{2n+1}_{eq} (S^{2n+1}, \alpha)\cong Ker(1-\alpha) | H_{2n+1} (S^{2n+1}).\]

Let $\tau $ be the fundamental class, then 
 \[[\tau] = \alpha ([\tau])\] 
  since $\alpha$ is an orientation preserving diffeomorphism. $\alpha$ induces identity map on $H_{2n+1}(S^{2n+1})$. Thus there holds the equality
   \[H^{2n+1}_{eq} (S^{2n+1},\alpha)=Ker(1-\alpha)|H_{2n+1}(S^{2n+1}) = \mathbb{C}.\]
$E_\infty ^{2n+1}(S^{2n+1}\rtimes_{\alpha_n} \mathbb{Z})$ contains at least one direct summand of $\mathbb{C}$ since
 \[E_\infty ^{2n+1}(S^{2n+1}\rtimes_{\alpha_n} \mathbb{Z}) = H^{2n+1} _{eq}(S^{2n+1},\alpha) \oplus H^{2n} _{coeq}(S^{2n+1},\alpha).\]

$\mathnormal{2)}$ When $k<0$, $H_{eq}^k (S^{2n+1},\alpha)$ and $H_{coeq}^k (S^{2n+1},\alpha)$ are $\{0\}$. Use the exact sequence again we have
\begin{eqnarray*}
 Coker\beta_0 &=&\{0\},\\
 Ker\gamma_0 &=&\{0\}, 
 \end{eqnarray*} 
i. e. 
\[H^{0}_{eq} (S^{2n+1},\alpha)\cong H_0 (S^{2n+1})/ {(1-\alpha) H_0 (S^{2n+1})} .\]
 
Let $dvol$ be the generator of $H^{2n+1}(S^{2n+1})$, i.e. the volume form. As $\alpha$ preserves orientation, the following equalities holds
  \[\int_{S^{2n+1}} \alpha(dvol) =  \int_{\alpha(S^{2n+1})} dvol =
  \int_{S^{2n+1}} dvol ,\]
 hence
  \[H_0 (S^{2n+1})/ {(1-\alpha)H_0 (S^{2n+1}) }=\mathbb{C}.\]
$E_\infty ^1(C^\infty (S^{2n+1}) \rtimes _ {\alpha_n}\mathbb{Z})$ contains at least one direct summand of $\mathbb{C}$.

From the following six-term exact sequence 
 \begin{displaymath}
 \xymatrix
 {HP^0 (C^\infty (S^{2n+1})) \ar[r] &HP^1(C^\infty(S^{2n+1}) \rtimes_{\alpha_n}\mathbb{Z} ) \ar[r]& HP^1(C^\infty(S^{2n+1})) \ar[d]^{1-\alpha} \\
  HP^0 (C^\infty (S^{2n+1})) \ar[u]^{1-\alpha}  & HP^0(C^\infty(S^{2n+1}) \rtimes_{\alpha_n}\mathbb{Z} )\ar[l] & HP^1(C^\infty(S^{2n+1})) \ar[l]}
 \end{displaymath}
 and the fact that $HP^0 (C^\infty (S^{2n+1}))= HP^1(C^\infty(S^{2n+1}))=\mathbb{C}$, 
 we know that 
 \[HP^0(C^\infty (S^{2n+1})\rtimes_{\alpha_n} \mathbb{Z}) = HP^1(C^\infty(S^{2n+1})\rtimes_{\alpha_n} \mathbb{Z}) =\mathbb{C} \oplus\mathbb{C}\]

Now it is obvious that 
  \begin{displaymath}
 E_\infty ^{2k+1}(S^{2n+1}\rtimes_{\alpha_n} \mathbb{Z}) = \left\{ \begin{array}{ll}
 \mathbb{C}, &k=0,n\\
 \{0\},         &else
\end{array}   \right.,                                                                                                                                                                                                                     
\end{displaymath}                                                                                                      
since
\[HP^1(C^\infty(S^{2n+1}) \rtimes_{\alpha_n}\mathbb{Z} ) = \oplus_{k=0}^{\infty} E_\infty ^{2k+1}(C^\infty (S^{2n+1}) \rtimes _ {\alpha_n}\mathbb{Z}).\]

    \end{proof}

Note that the proof above is inspired by Nest's work in \cite{nest}. 

\section{Another example}	

People may think that smooth crossed products of diffeomorphisms depend on manifolds only. This is not true. Our observation is based on the beautiful result of Brenken-Cuntz-Elliott-Nest \cite{bcen}. They utilized the theory of cyclic cohomology to obtain their result. 

Let us recall some notions about noncommutative tori of three dimensions. Let $\eta: \mathbb{Z}^3\wedge \mathbb{Z}^3 \to \mathbb{T}$ be an antisymmetric bicharachter. Denote $\eta(e_i\wedge e_j), i,j=1,2,3$ as $\eta_{i,j}$, where $e_i,e_j$ are the geneorators of $\mathbb{Z}^3$. $A_\eta$ is the $C^*$-algebras generated by three unitaries $u_1,u_2,u_3$ under the relation
\[u_i u_j = \eta_{i,j} u_j u_i.\]
Call $\eta$ nondegenerate if and only if $\eta(\mathbb{Z}^3 \wedge g )=1, g\in \mathbb{Z}^3$ implies $g=0$.
Let 
$\mathscr{S}(\mathbb{Z}^3)$ be the space of rapidly decreasing sequences on $\mathbb{Z}^3_{u_1,u_2,u_3}$. Endow it with the topology given by seminorms 
\[\|(x_a)_{a\in \mathbb{Z}^3_{u_1,u_2,u_3}}\|_k = \sup_{a\in \mathbb{Z}^3_{u_1,u_2,u_3}} (1+|a|^k)|x_a|.\]
 $\mathcal{A}_\eta$ is the completion of $\mathscr{S}(\mathbb{Z}^3)$ under this topology.

\begin{theo}[\cite{bcen}]
	Let $\eta$ and $\eta'$ be antisymmetric bicharacters on 	$\mathbb{Z}^3$, both nondegererate. Then $\mathcal{A}_\eta \cong \mathcal{A}_{\eta'}$  if and only if $\eta \cong \eta'$.
\end{theo}
We will consider a special class of $\mathcal{A}_{\eta}$ ( $A_{\eta}$ ) which can be described as smooth crossed product ( $C^*$ crossed product) algebras of diffeomorphisms of $\mathbb{T}^2$ by $\mathbb{Z}$. Let $\beta$ be a diffeomorphism of $\mathbb{T}^2$ such that $\beta (x_1, x_2) = (e^{2\pi i\theta_1}x_1,e^{2\pi i\theta_2}x_2)$, where $\theta_1$ and $\theta_2$ are two rationally independent irrational numbers. We choose $X_1 =\frac{\partial}{\partial x_1} $ and $X_2 = \frac{\partial}{\partial x_2}$  to obtain the seminorms on $C^\infty(\mathbb{T}^2)$. Notice that $\rho_k(n)\equiv(1+2|n|)^k$ since $\|\beta^t\|_n= 1$.

By definitions and basic Fourier analysis we know that
\begin{eqnarray*} 
C^\infty(\mathbb{T}^2)\rtimes_\beta \mathbb{Z}&\cong & \mathcal{A}_\eta,\\
C(\mathbb{T}^2)\rtimes_\beta \mathbb{Z} &\cong & A_\eta
\end{eqnarray*}
where  
\begin{displaymath}
	\mathbf{\eta}=
	\left(\begin{array}{ccc}
	1&1&	e^{-2\pi i\theta_1}\\
	1&1&    e^{-2\pi i\theta_2}\\
	e^{2\pi i\theta_1} & e^{2\pi i\theta_2} &1 
	\end{array}
	\right).
\end{displaymath}

Now apply theorem 4.1 we have 
\begin{exam}
$\beta$ and $\beta'$ are two minimal unique ergodic diffeomorphism of $\mathbb{T}^2$ such that:
 \begin{eqnarray*}
 \beta (x_1, x_2) &=& (e^{2\pi i\theta_1}x_1,e^{2\pi i\theta_2}x_2),\\
 \beta' (x_1, x_2) &=& (e^{2\pi i\theta_1'}x_1,e^{2\pi i\theta_2'}x_2).
 \end{eqnarray*}
 $\theta_1,\theta_2,\theta_1',\theta_2'	$ are rationally independent irrational numbers. Since their corresponding bicharacters are not isomorphic to each other and are both nondegenerated, $\beta$ and $\beta'$ give different smooth crossed product algebras.
\end{exam}

\begin{rem}
Note that one can also find a couple of diffeomorphisms of $S^1$ giving different smooth crossed product algebras in essentially the same way. 	
\end{rem}

\end{document}